\documentstyle[11pt]{article}

\newtheorem{thm}{Theorem}[section]
\newtheorem{lem}[thm]{Lemma}
\newtheorem{cor}[thm]{Corollary}
\newtheorem{pro}[thm]{Proposition}
\newtheorem{ex}[thm]{Example}
\newtheorem{rmk}[thm]{Remark}
\newtheorem{defi}[thm]{Definition}

\newcommand{\lon }{\longrightarrow }
\newcommand{\be }{\begin{eqnarray*}}
\newcommand{\ee }{\end{eqnarray*}}

\input amssym

\def\nin{\noindent}

\newcommand{\pf}{\noindent{\bf Proof.}\ }
\newcommand{\qed}{\begin{flushright} $\Box$\ \ \ \ \ \
                  \end{flushright}}

\newcommand{\reals}{{\Bbb R}}

\newcommand{\frakg}{{\frak g}}

\newcommand{\cinf}{C^{\infty}}

\newcommand{\smalcirc}{\mbox{\tiny{$\circ $}}}

\def\description label#1{\hfil\bf[#1]\hfil}
\newcommand{\R}{\reals}

\newcommand{\Ga}{\Gamma}

\def\Clin{C^\infty_{\ell in}}

\def\mvf{multiplicative vector field}

\def\cdo{covariant differential operator}

\def\CDO{\mathop{\rm CDO}}
\def\im{\mathop{\rm im}}
\def\End{\mathop{\rm End}}
\def\fms#1{\Omega^1(#1)}

\def\ST{\ \vert\ }

\def\sdp{\mathbin{\hbox{$\mapstochar\kern-.3333em\times$}}}
\def\pds{\mathbin{\hbox{$\times\kern-.55em\mapstochar\,$}}}

\newcommand{\wed}{\mathbin{\lower1.5pt\hbox{$\scriptstyle{\wedge}$}}}

\let\Tilde=\widetilde
\let\Bar=\overline
\let\Vec=\overrightarrow
\let\ceV=\overleftarrow
\let\Hat=\widehat
\let\isom=\cong
\let\sol=\bullet
\def\chigh{{\raise1.5pt\hbox{$\chi$}}}
\let\phi=\varphi
\def\tilalpha{\widetilde\alpha}
\def\tilbeta{\skew6\widetilde\beta}
\def\tilone{\widetilde 1}
\def\til0{\Tilde{0}}

\def\dpl{\mathbin{+\hskip-6pt +\hskip4pt}}
\def\dminus{\raise2pt\hbox{\vrule height1pt width 2ex}\hskip3pt}
\def\dtimes{\mathbin{\hbox{\huge.}}}

\def\llangle{\langle\!\langle}
\def\rrangle{\rangle\!\rangle}

\def\pback#1{\mathbin{{{\lower1.2ex\hbox{$\times$}}\atop #1}}}

\def\ddt#1{\left.\frac{d}{dt}#1\right|_0}

\def\upa{\uparrow}

\def\hatt{\Hat{\phantom{X}}}

\def\vlra{\hbox{$\,-\!\!\!-\!\!\!-\!\!\!-\!\!\!-\!\!\!
-\!\!\!-\!\!\!-\!\!\!-\!\!\!-\!\!\!\longrightarrow\,$}}

\def\lrah{\hbox{$\,-\!\!\!-\!\!\!
-\!\!\!-\!\!\!-\!\!\!-\!\!\!-\!\!\!\longrightarrow\,$}}

\def\surj{-\!\!\!-\!\!\!-\!\!\!\gg}

\def\inj{>\!\!\!-\!\!\!-\!\!\!-\!\!\!>}

\def\gpd{\,\lower1pt\hbox{$\longrightarrow$}\hskip-.24in\raise2pt
             \hbox{$\longrightarrow$}\,}

\def\lgpd{\,\lower1pt\hbox{$\vlra$}\hskip-1.02in\raise2pt\hbox{$\vlra$}\,}

\def\llgpd{\,\lower1pt\hbox{$\vvlra$}\hskip-1.3in\raise2pt\hbox{$\vvlra$}\,}

\def\kwo#1{\lq\lq #1\rq\rq}

\begin{document}

\title{{\bf Classical lifting processes and multiplicative vector fields}
\thanks{1991 {\em Mathematics
Subject Classification.} Primary 58F05. Secondary 17B66, 22A22, 53C99,
58H05.}}

\author{KIRILL C. H. MACKENZIE\\
        School of Mathematics and Statistics\\
        University of Sheffield\\
        Sheffield, S3 7RH, England\\
        {\sf email: K.Mackenzie@sheffield.ac.uk}\\
        and \\
        PING XU
        \thanks{Research partially supported by NSF grants DMS92-03398
        and DMS95-04913, and a NSF postdoctoral fellowship.}\\
   Department of Mathematics\\
   The Pennsylvania State University \\
   University Park, PA 16802, USA\\
   {\sf email: ping@math.psu.edu}}

\date{February 1, 1997}

\maketitle

\begin{abstract}
We extend the calculus of multiplicative vector fields and differential
forms and their intrinsic derivatives from Lie groups to Lie groupoids;
this generalization turns out to include also the classical process of
complete lifting from arbitrary manifolds to tangent and cotangent bundles.
Using this calculus we give a new description of the Lie bialgebroid
structure associated with a Poisson groupoid.
\end{abstract}

\section{Introduction}

Multiplicative multivector fields and forms play an important technical
role in Poisson group theory \cite{LuW:1990}, \cite{Weinstein:1990}:
most basically, the Poisson tensor is itself multiplicative and its
intrinsic derivative gives the bracket on the Lie algebra dual. If one
extends this approach to Poisson groupoids---an extension we did not
contemplate in \cite{MackenzieX:1994}, \cite{MackenzieX2}---it turns out to also
generalize the classical concept of complete lift.

Complete and vertical lifts of vector fields and differential forms from a
manifold to its tangent and cotangent bundles are dealt with comprehensively
in the work of Yano and Ishihara \cite{YanoI}. The existence of the complete
lifting processes is one of the fundamental features which set the geometry
of tangent and cotangent bundles apart from that of general vector and
principal bundles, and it is consequently reasonable to expect that it
extend to Lie algebroids and their duals. In fact it extends only when
the Lie algebroid is integrable: the complete lift of a vector field on a
manifold $M$ to $TM$ is effectively the intrinsic derivative of an
associated multiplicative vector field on the pair groupoid
$M\times M$. In general, a multiplicative vector field on a groupoid $G$
induces vector fields on its Lie algebroid $AG$ and its dual $A^*G$. We
characterize the former in terms of the Lie algebroid structure of
$TAG\lon TM$, and the latter in terms of the Poisson structure on $A^*G$
(Theorem \ref{thm:movf}). Multiplicative forms on $G$ likewise induce forms
on $AG$ and $A^*G$; in this paper, where we are concerned with applications
to Poisson groupoids, we deal only with 1-forms.

The importance of these processes is that, by slightly generalizing the
notion of multiplicative vector field, and incorporating the vertical lifts,
we obtain a complete description, well adapted to the geometry of the
groupoid, of the vector fields on the Lie algebroid of a Lie groupoid.
In particular, we give in \S\ref{sect:fafopg} a new description of the
Poisson structure on the Lie algebroid of a Poisson groupoid, and thus of
the Lie algebroid structure on the dual.

We begin in \S\ref{sect:lvf} with the case of vector fields on a vector
bundle $(A,q,M)$. Here it is natural to consider those vector fields
$X\colon A\to TA$ which are vector bundle morphisms with respect to the
vector bundle structure of $TA$ over $TM$. Such {\em linear vector fields}
correspond to covariant differential operators
$D\colon \Gamma A\to \Gamma A$ in the same way that, given a connection
$\nabla$ in $A$, the horizontal lift of a vector field $x$ on $M$
corresponds to $\nabla_x\colon\Gamma A\to\Gamma A$. Secondly, any section
$X$ of $A$ induces a vector field $X^\upa$ on $A$: if $A = TM$ then $X^\upa$
is the vertical lift of $X$; in the general case we call it the {\em core
vector field} corresponding to $X$. The core vector fields and the linear
vector fields together generate ${\cal X}(A)$. We believe that the material
of this section is mostly folklore.

In \S\ref{sect:vfolg} we give a similar calculus on general Lie groupoids.
Given a multiplicative vector field $\xi$ on a groupoid $G\gpd M$,
application of the Lie functor leads to a linear vector field $\Tilde{\xi}$
on $AG$ for which the corresponding differential operator $D_\xi$ is a
derivation of the Lie algebroid bracket. In the case where $G$ is a pair
groupoid $M\times M$, the process $\xi\mapsto\Tilde{\xi}$ is the complete
lifting of \cite{YanoI}. For a Lie group $G$ it is linearization at the
identity. In order to obtain a complete description of the vector fields
on $AG$ we weaken the multiplicativity condition to what we call a {\em
star vector field}; for any Lie groupoid $G$ the $\Tilde{\xi}$, for $\xi$
a star vector field, and the core vector fields $X^\upa$, for
$X\in\Gamma AG$, generate ${\cal X}(AG)$.

For any vector bundle $A$, there is a bijective correspondence between
the linear vector fields on $A$ and the linear vector fields on the dual
$A^*$. Thus for a Lie algebroid dual $A^*G$ we obtain in \S\ref{sect:vfolad}
a description of the vector fields on $A^*G$ in terms of the \kwo{duals}
of the $\Tilde{\xi}$, for $\xi$ a star vector field on $G$, and the core
vector fields $\phi^\upa$ for $\phi\in\Gamma A^*G$. More generally, we show
in Theorem~\ref{thm:movf} that for any abstract Lie algebroid $A$, a linear
vector field $\xi$ is a Lie algebroid morphism $A\to TA$ if and only if
the corresponding vector field on $A^*$ is Poisson (with respect to the
dual Poisson structure on $A^*$), and this is so if and only if $D_\xi$
is a derivation of the Lie algebroid bracket.

In \S\ref{sect:flg} and \S\ref{sect:folad} we give a comparable analysis
for differential 1--forms on Lie algebroids and their duals.

For a Poisson groupoid $G\gpd P$, the 1--forms have a bracket structure
reflecting the fact that $T^*G\to G$ is a Lie algebroid. In the final
\S\ref{sect:fafopg} we show that 1--forms on a Poisson groupoid admit a
calculus similar to that of \S\ref{sect:vfolg}, with the differential
operators on $AG$ now replaced by operators $\fms{P}\to\fms{P}$; for
example, if $\Phi\in\fms{G}$ is multiplicative, then $D_\Phi$ is a
derivation of the Poisson bracket on $\fms{P}$. In these terms we obtain
in Theorem~\ref{thm:last}  a complete description of the bracket structure
on $\fms{G}$.

It is worth noting that the treatment we give here is entirely coordinate
free, and so may offer something new even in the classical case. In the
early stages of the work we were unaware of \cite{YanoI} and our approach
was chiefly influenced by the few brief remarks in \cite{Dieudonne:IV} and
\cite{AbrahamM}.

Some of the material of this paper was announced in \cite{Mackenzie:AMI}.

We are grateful to Alan Weinstein for some helpful comments, and to the
Isaac Newton Institute, where much of this work was done, for its excellent
facilities and environment.

\section{Linear vector fields}  \label{sect:lvf}

We first need some preliminaries concerning linear vector fields on vector
bundles. These are, in effect, homogeneous vector fields of degree 1 in
the sense of \cite{Besse}, but the point of view we adopt here will be
important in what follows. Consider a fixed vector bundle $(A,q,M)$. Recall
from \cite{MackenzieX:1994} and references given there the {\em tangent
double vector bundle}
\begin{equation}                         \label{tdvb}
\matrix{&&T(q)&&\cr
        &TA&\vlra&TM&\cr
        &&&&\cr
     p_A&\Bigg\downarrow&&\Bigg\downarrow&p\cr
        &&&&\cr
        &A&\vlra&M.&\cr
        &&q&&\cr}
\end{equation}
We recall the notation of \cite{MackenzieX:1994}. In $A$ and $TM$, we use
standard notation, with the zero of $A$ over $m\in M$ being $0^A_m$, and the
zero of $TM$ over $m$ being $0_m^T$.

We denote elements of $TA$ by $\xi,\eta,\zeta\ldots,$ and we write
$(\xi;X,x;m)$ to indicate that $X = p_A(\xi),\ x = T(q)(\xi),$ and
$m = p(T(q)(\xi)) = q(p_A(\xi))$. With respect to the
tangent bundle structure $(TA, p_A, A)$, we use standard notation: $+$ for
addition, $-$ for subtraction and juxtaposition for scalar multiplication. The
notation $T_X(A)$ will always denote the fibre $p_A^{-1}(X)$, for
$X\in A$, with respect to this bundle. The zero element in $T_X(A)$ is denoted
$\Tilde{0}_X$. We refer to this bundle structure as the {\em $p_A$-bundle
structure}.

With respect to the {\em $T(q)$-bundle structure}, $(TA,T(q),TM)$, we use
$\dpl$ for addition, $\dminus$ for subtraction, and $\dtimes$ for scalar
multiplication. This addition and scalar multiplication on $TA$ are precisely
the tangents of the addition and scalar multiplication in $A$.
The fibre over $x\in TM$ will always be denoted
$T(q)^{-1}(x)$, and the zero element of this fibre is $T(0)(x)$.

For each $m\in M$, the tangent space $T_{0_m}(A_m)$ identifies canonically
with $A_m$; we denote the element of $T_{0_m}(A_m)$ corresponding to $X\in
A_m$ by $\Bar{X}$. The elements $\Bar{X},\ X\in A$, form the {\em core} of
$TA$. Note that, for $X,Y\in A_m$ and $t\in\R$,
$$
\Bar{X} + \Bar{Y} = \Bar{X + Y} = \Bar{X} \dpl \Bar{Y},\qquad
t\Bar{X} = \Bar{tX} = t\dtimes\Bar{X}.
$$

Given a morphism of vector bundles $\phi\colon A'\to A,\ f\colon M'\to M$,
we denote the pullback of $A$ across $f$ by $f^!A$, and the induced morphism
$A'\to f^!A$ over $M'$ by $\phi^!$. Associated to the double vector bundle
structure on $TA$ are the two core
sequences:
\begin{equation}           \label{eq:pAcore}
q^!A\stackrel{\tau}{\inj} TA \stackrel{T(q)^!}{\surj}q^!TM,\qquad\qquad
p^!A\stackrel{\upsilon}{\inj} TA \stackrel{p_A^!}{\surj} p^!A,
\end{equation}
over $A$ and $TM$ respectively. Here $\tau$ and $\upsilon$ are the maps
$$
\tau(X,Y) = \Tilde{0}_X \dpl \Bar{Y},\qquad
\upsilon(x,Y) = T(0)(x) + \Bar{Y};
$$
we call $\tau$ and $\upsilon$ the {\em translation maps.} $\tau$ assigns to
$(X,Y)\in A_m\times A_m$ the element of $T_X(A_m)$ which has its tail at $X$,
and is parallel to $Y$.

Given $X\in\Gamma(A)$ define a vector field $X^\upa$ on $A$ by
$X^\upa(Y) = \tau(Y,X(qY)),\ Y\in A$. (In \cite{MackenzieX:1994} we used
the notation $\breve{X}$.) Then
$$
X^\upa(F)(Y) = \ddt{F(Y+tX(qY))}
$$
for $F\in \cinf(A),\ Y\in A$, and so
\begin{equation}                         \label{eq:cvf}
X^\upa(f\circ q) = 0,\qquad
X^\upa(l_\phi) = \langle\phi,X\rangle\circ q,\qquad
[X^\upa,Y^\upa] = 0,
\end{equation}
for $f\in \cinf(M),\ X,Y\in\Gamma(A),\ \phi\in\Gamma(A^*)$. Here $\ell_\phi\in
\cinf(A)$ is the fibrewise linear function determined by $\phi$, namely
$X\mapsto\langle\phi(qX),X\rangle.$ Note also that
$(fX)^\upa = (f\circ q)X^\upa$. We call $X^\upa$ the {\em core vector field}
corresponding to $X\in\Ga A$; if $A = TM$, it is the vertical lift of $X$.

A section $X\in\Gamma(A)$ also induces a section $\Hat{X}$ of $T(q)$ by
$\Hat{X}(x) = \upsilon(x,X(px))$. Note that
$$
(X + Y)\hatt = \Hat{X}\dpl\Hat{Y},\qquad
(fX)\hatt = (f\circ p)\dtimes\Hat{X},
$$
for $X,Y\in\Gamma(A), f\in \cinf(M)$.

\begin{defi}
A {\em linear vector field on $A$} is a pair $(\xi, x)$, where $\xi$
is a vector field on $A$, and $x$ is a vector field on $M$, such that
\begin{equation}                         \label{lvf}
\matrix{&&\xi&&\cr
        &A&\vlra&TA&\cr
        &&&&\cr
       q&\Bigg\downarrow&&\Bigg\downarrow&T(q)\cr
        &&&&\cr
        &M&\vlra&TM,&\cr
        &&x&&\cr}
\end{equation}
is a morphism of vector bundles.
\end{defi}

In particular, $\xi$ projects under $q$ to $x$, and
$$
\xi(X+Y) = \xi(X)\dpl \xi(Y),\qquad \xi(tX) =
t\dtimes\xi(X),
$$
for all $X,Y\in A,\ t\in\R$. In what follows, $\Clin(A)$ denotes the
subset of $\cinf(A)$ consisting of the fibrewise linear functions
$\ell_\phi\colon A\to\R,\ \phi\in\Ga A^*$.

Given a connection in $A$, the horizontal lift of any vector field $x$ on $M$
is a linear vector field on $A$ over $x$. Core vector fields
$X^\upa,\ X\neq 0,$ are not linear.

\begin{pro}                         \label{TFAE(lvf)}
Let $\xi$ be a vector field on $A$. The following are equivalent:
\begin{enumerate}
\item $(\xi,x)$ is a linear vector field on $A$ for some vector field $x$ on
$M$;                                             \label{vf}
\item $\xi\colon \cinf(A)\to \cinf(A)$ sends     \label{functs}
$\Clin(A)$ into $\Clin(A)$ and sends $q^*\cinf(M)$ into $q^*\cinf(M)$;
\item $\xi$ has flows $\phi_t$ which are (local) vector bundle morphisms
$A\to A$ over a flow $f_t$ on $M$.               \label{flow}
\end{enumerate}
\end{pro}

\nin {\sc Proof:}\
(i)$\Longrightarrow$(ii). Suppose first that $\xi, \eta\in T(A)$
have $T(q)(\xi) = T(q)(\eta)$, and that $F\in \Clin(A)$. Write
$\xi = \left.\frac{d}{dt}X_t\right|_0,\ \eta = \left.
\frac{d}{dt}Y_t\right|_0$, where $q(X_t) = q(Y_t)$ for $t$ close to $0\in\R$.
Then
\begin{eqnarray*}
(\xi\dpl\eta)(F)(X_0+Y_0)
   & = & \left.\frac{d}{dt}F(X_t + Y_t)\right|_0\\
   & = & \left.\frac{d}{dt}F(X_t)\right|_0 +
                      \left.\frac{d}{dt}F(Y_t)\right|_0\\
   & = & \xi(F)(X_0) + \eta(F)(Y_0).
\end{eqnarray*}
>From this, and the similar result for scalar multiplication, it
follows that a linear vector field $\xi$ maps $\Clin(A)$ into
$\Clin(A)$. Since $\xi$ is projectable under $q$, it is clear that
$\xi$ maps $q^*\cinf(M)$ into $q^*\cinf(M)$.

(ii)$\Longrightarrow$(iii).
>From the assumption that $\xi$ sends $q^*\cinf(M)$ into itself, and a standard
result, it follows that $\xi$ projects to a vector field $x$ on $M$, and
hence $q\circ\phi_t = f_t\circ q$ for all $t$, where $\phi_t$ and $f_t$ are
flows for $\xi$ and $x$ respectively. It remains to prove that the $\phi_t$
are linear and since this is a local question, we may assume that
$A = M\times V$ and that $\xi$ and $x$ have global flows. The details are
straightforward.

(iii)$\Longrightarrow$(i). Suppose, for simplicity, that $\xi$ has
a global flow $\phi_t$ by vector bundle morphisms over a global flow $f_t$ on
$M$. Then $\xi$ certainly projects under $q$ to the vector field $x$
generated by $f_t$. For $X,Y\in A$ with $q(X) = q(Y)$,
\begin{eqnarray*}
\xi(X+Y) & = & T(t\mapsto \phi_t(X+Y))_0(1)\\
             & = & T(t\mapsto \phi_t(X) + \phi_t(Y))_0(1)\\
             & = & T(+)(T(t\mapsto\phi_t(X))_0(1),
                                    T(t\mapsto\phi_t(Y))_0(1))\\
             & = & T(+)(\xi(X),\xi(Y))\\
             & = & \xi(X)\dpl\xi(Y),
\end{eqnarray*}
and similarly one proves that $\xi(t X) = t\dtimes\xi(X)$
for $t\in\R$.
\qed

\begin{cor}                 \label{cor:tfae}
If $(\xi,x)$ and $(\eta,y)$ are linear vector fields, then
$([\xi,\eta],[x,y])$ is also.
\end{cor}

Now consider a linear vector field $(\xi,x)$. Since $\Clin(A)$ is
canonically isomorphic (as a $\cinf(M)$ module) to $\Ga A^*$, condition (ii)
above shows that $\xi$ induces a map
\begin{equation}                            \label{eq:D*}
D^{(*)}_\xi\colon\Ga A^*\to\Ga A^*,\quad \mbox{such that}\quad
\ell_{D^{(*)}_\xi(\phi)} = \xi(\ell_\phi).
\end{equation}
Clearly, $D^{(*)}_\xi$ is additive, and it is easily checked that
$D^{(*)}_\xi(f\phi) = fD^{(*)}_\xi(\phi) + x(f)\phi$ for
$f\in \cinf(M)$. Thus $D^{(*)}_\xi$ is a covariant differential operator
on $A^*$ (\cite[III\S2]{Mackenzie:LGLADG}), and we have a map, which is linear
over $\cinf(M)$ and bracket-preserving,
\begin{equation}              \label{isom}
(\xi,x)\mapsto D^{(*)}_\xi,\quad \Ga^{LIN}TA\to \Ga\CDO(A^*),
\end{equation}
where $\Ga^{LIN}TA$ is the set of linear vector fields on $A$ and $\CDO(A^*)$
is the vector bundle on $M$ whose sections are \R-linear maps
$D\colon\Ga A^*\to \Ga A^*$ such that $D(f\phi) = fD(\phi) + x(f)\phi$ for all
$f\in \cinf(M),\ \phi\in\Ga A^*,$ and some fixed $x\in{\cal X}(M)$. With
respect to the commutator bracket, $\CDO(A^*)$ is a Lie algebroid on $M$ with
anchor $D\mapsto x$.

By Corollary \ref{cor:tfae}, $\Ga^{LIN}TA$ is closed under the bracket on
${\cal X}(A)$. Let $a$ denote the map $(\xi,x)\mapsto x$. Then the kernel of
$a$ consists of those vertical vector fields on $A$ which are linear. Any
vertical vector field $\xi$ can be written as $\xi(X) = \tau(X,{\sf X}(X))$,
where ${\sf X}\colon A\to A$ has $q\circ{\sf X} = {\sf X}$, and $\xi$ is
linear iff ${\sf X}$ is a vector bundle morphism. Thus the kernel of $a$
can be identified with $\Ga\End(A)$. Since this is the module of sections of
a vector bundle over $M$, it follows that $\Ga^{LIN}TA$ is also the module of
sections of a vector bundle on $M$; we denote this by $T^{LIN}A$. We can now
regard (\ref{isom}) as a morphism of Lie algebroids from $T^{LIN}A$ to
$\CDO(A^*)$.

\begin{pro}                                \label{pro:isom}
The morphism $T^{LIN}A\to\CDO(A^*)$ just defined is an isomorphism of Lie
algebroids.
\end{pro}

\nin {\sc Proof:}\
It suffices to prove (see \cite[III~2.8]{Mackenzie:LGLADG}) that the
restriction of this map to the kernels of the anchor maps (the adjoint
bundles) is an isomorphism, and this restriction is the canonical
identification of $\End(A)$ with $\End(A^*)$.
\qed

Now each covariant differential operator $D^{(*)}$ on $A^*$ corresponds to a
covariant differential operator $D$ on $A$ by
$$
\langle D^{(*)}(\phi),X\rangle = a(D)(\langle\phi,X\rangle) -
\langle\phi,D(X)\rangle,
$$
where $\phi\in\Ga A^*,\ X\in\Ga A$. This defines an
isomorphism of Lie algebroids $\CDO(A)\isom\CDO(A^*)$. In sum, we have
obtained a correspondence between linear vector fields on $A$, and covariant
differential operators on either $A$ or $A^*$. Letting $D_\xi$ denote the
element of $\Ga\CDO(A)$ corresponding to $D^{(*)}_\xi$, we have
\begin{equation}                         \label{eq:D}
\langle\phi,D_\xi(X)\rangle = x(\langle\phi,X\rangle) -
      \xi(\ell_\phi)\circ X.
\end{equation}

\begin{pro}                            \label{prop:HigginsM:1990a}
Let $(\xi,x)$ be a linear vector field on $A$ and let $D = D_\xi$ be the
corresponding element of $\Ga\CDO(A)$. Then, for all $X\in\Ga A$ and
$m\in M$,
$$
\tau(X(m),D(X)(m)) = T(X)(x(m)) - \xi(X(m)).
$$
If $\phi_t$ is a flow for $\xi$ near $X(m)$ and $f_t$ the corresponding flow
for $x$ near $m$, then
$$
\upsilon(x(m),D(X)(m)) =
\left.\frac{d}{dt}\left(X(f_t(m)) - \phi_t(X(m))\right)\right|_0.
$$
\end{pro}

\pf
Let $Y = D(X)(m)$ and let $Z$ be the RHS of the first equation. It suffices
to prove that $\tau(X(m),Y) = \tau(X(m),Z)$, both of which are vertical
tangent vectors to $A$ at $X(m)$. The functions on $A$ are generated by
those of the form $\ell_\phi$ for $\phi\in\Ga A^*$ and those of the form
$f\circ q$ for $f\in \cinf(M)$. Since both vectors are vertical, they coincide
on all $f\circ q$. Now, for any $\phi\in\Ga A^*$,
$$
\tau(X(m),Y)(\ell_\phi) = \langle\phi,Y\rangle
= x(m)(\langle\phi,X\rangle) - \xi(\ell_\phi)(X(m))
$$
and
$$
\tau(X(m),Z)(\ell_\phi) = T(X)(x(m))(\ell_\phi) -
\xi(X(m))(\ell_\phi),
$$
whence the result.

For the second equation, note first that $T(X)(x(m))$ and $\xi(X(m))$
have the same two projections, and therefore
$$
\upsilon(x(m),D(X)(m)) = T(X)(x(m))\dminus\xi(X(m)).
$$
>From this the second equation follows.
\qed

If, in Proposition \ref{prop:HigginsM:1990a}, $X(m) = 0$ for a specific $m\in M$, then
$D(X)(m)$ depends only on $x,X$ and $m$ and can be identified with
$T(X)(x(m))$; the map $T_m(M)\to A_m,\ x\mapsto T(X)(x),$ is the
{\em intrinsic derivative} of \cite{GolubitskyG}. If, in addition, $A = TM$,
this map is the {\em linearization of a vector field at a singularity}
in the sense of \cite[p.72]{AbrahamM}.

Notice that there is now a bijective correspondence between linear vector
fields on $A$ and linear vector fields on $A^*$. It follows from
Proposition \ref{xiH} below that if
$(\xi,x)\in{\cal X}^{LIN}(A)$ and $(\xi_*,x)\in{\cal X}^{LIN}(A^*)$
correspond in this way, and if $\phi_t$ is a (local) flow for $\xi$, then
$\phi_{-t}^*$ is a (local) flow for $\xi_*$.

Next recall the {\em tangent pairing} between $TA\to TM$ and $T(A^*)\to TM$
of \cite[5.3]{MackenzieX:1994}: given ${\goth X}\in T(A^*)$ and $\xi\in TA$
with $T(q)({\goth X}) = T(q_*)(\xi)$ we can write
${\goth X} = \ddt{\phi_t}\in T(A^*)$ and $\xi = \ddt{X_t}\in TA$
where $X_t\in A$ and $\phi_t\in A^*$ have $q_*(\phi_t) = q(X_t)$ for $t$
near zero. Now the tangent pairing $\llangle\ ,\ \rrangle$ is defined by
\begin{equation}                             \label{eq:pairing}
\llangle{\goth X},\xi\rrangle = \ddt{\langle\phi_t,X_t\rangle}.
\end{equation}

\begin{pro}                    \label{prop:llangle}
{\rm \cite[6.3]{MackenzieX:1994}}
Given $(\xi;X_m,x;m)\in TA$ and $({\goth X};\phi_m,x;m)\in T(A^*)$, let
$X\in\Gamma(A)$ and $\phi\in\Gamma(A^*)$ be any sections taking the values
$X_m$ and $\phi_m$ at $m$. Then
$$
\llangle{\goth X},\xi\rrangle = {\goth X}(l_X) + \xi(l_\phi)
                                  - x(\langle\phi,X\rangle).
$$
\end{pro}

>From this the following characterization of the dual of a linear vector field
follows easily.

\begin{pro}                    \label{xiH}
Let $(\xi,x)$ be a linear vector field on $A$. For $\phi\in A^*_m,\ m\in M$,
the value $H = \xi_*(\phi)$ of the corresponding linear vector field on
$A^*$ is the unique element of the form $(H;\phi,x(m);m)\in T(A^*)$ for which
$$
\llangle H,\xi(X)\rrangle = 0
$$
for all $X\in A_m$.
\end{pro}

\section{Vector fields on Lie groupoids}         \label{sect:vfolg}

Now consider a Lie groupoid $G$ on base $M$, and its Lie algebroid $AG$;
see \cite{Mackenzie:LGLADG} for the basic facts concerning Lie groupoids and Lie
algebroids. Applying the tangent functor to the operations in $G\gpd M$
yields the {\em tangent groupoid} $TG\gpd TM$. If $\kappa\colon G*G\to G$ is
the multiplication, then the multiplication in $TG$ is
$X\sol Y = T(\kappa)(X,Y)$.

\begin{pro} {\rm \cite[2.6]{Xu:1995}}            \label{pro:sol}
Let $X\in T_g(G)$ and $Y\in T_h(G)$ have $T(\alpha)(X) = T(\beta)(Y)  = x$.
Then
$$
X\sol Y = T(L_{{\cal X}})(Y) + T(R_{{\cal Y}})(X) -
T(L_{{\cal X}})T(R_{{\cal Y}})(T(1)(x))
$$
where ${\cal X,Y}$ are any bisections of $G$ with ${\cal X}(\alpha g) = g,\
{\cal Y}(\alpha h) = h$.
\end{pro}

Here 1 is the identity map embedding the base manifold $M$ into $G$.
A {\em bisection} \cite{CDW} of $G$ is a submanifold ${\cal X}$ of $G$ such
that $\alpha\colon{\cal X}\to M$ and $\beta\colon{\cal X}\to M$ are both
diffeomorphisms. Bisections are naturally identified with admissible sections
\cite[II\S5]{Mackenzie:LGLADG} and we will use both formulations without
comment.
(If $G$ is locally trivial, the bisections are the gauge transformations
of the corresponding principal bundle.)

The tangent bundle projection $p_G\colon TG\to G$ is a groupoid morphism over
$p\colon TM\to M$ and applying the Lie functor yields a morphism of Lie
algebroids
\begin{equation}                       \label{diag:LATG}
\matrix{&&q_{TG}&&\cr
        &ATG&\lrah&TM&\cr
        &&&&\cr
  A(p_G)&\Big\downarrow&&\Big\downarrow&p_M\cr
        &&&&\cr
        &AG&\lrah&M,&\cr
        &&q_G&&\cr}
\end{equation}
where $q_{TG}$ and $q_G$ are the bundle projections of the Lie algebroids.
In fact $ATG\to AG$ can be given a vector bundle structure, by applying the
Lie functor to the operations in $TG\to G$ \cite[\S7]{MackenzieX:1994}. This
gives (\ref{diag:LATG}) the structure of a double vector bundle with core
$AG$.

\begin{thm}                        \label{thm:j}
{\rm \cite[7.1]{MackenzieX:1994}}
Let $G$ be a Lie groupoid on base $M$. Then there is a canonical isomorphism
of double vector bundles $j_G\colon TAG\to ATG$, where $ATG$ is as above and
$TAG$ is the tangent double vector bundle of $AG\to M$, which induces the
identities on the side bundles $AG$ and $TM$ and on the cores $AG$.
\end{thm}

\begin{defi}                               \label{defi:mvf}
A {\em multiplicative vector field} on $G$ is a pair of vector fields $(\xi,x)$
where $\xi\in{\cal X}(G),\ x\in{\cal X}(M)$, such that $\xi\colon G\to TG$ is
a morphism of groupoids over $x\colon M\to TM$.
\end{defi}

For a group $G$, we must of course have $x = 0$ and so \ref{defi:mvf}
coincides with the definition in \cite[1.3]{LuW:1990}. For a pair groupoid
$M\times M$, the multiplicative vector fields are those of the form
$x\times x$, where $x\in{\cal X}(M)$.

\begin{ex}\rm
Given $X\in\Ga AG$, the right- and left-invariant vector fields corresponding
to $X$ are defined by $\Vec{X}(g) = T(R_g)(X(\beta g))$ and $\ceV{X}(g) =
T(L_g)T(i)(X(\alpha g))$, where $i\colon G\lon G$ is the inversion in $G$.
Write $\xi = \Vec{X} + \ceV{X}$. Then $\xi$ is a multiplicative vector field
over $a(X)$, the anchor of $X$.
\end{ex}

\begin{pro}                       \label{pro:mvf}
Let $\xi$ be a vector field on a Lie groupoid $G$ and $x$ a vector field on
its base $M$. Then the following are equivalent:
\begin{enumerate}
\item $(\xi,x)$ is a multiplicative vector field on $G$;

\item The flows $\phi_t$ of $\xi$ are (local) Lie groupoid automorphisms
over the flows $f_t$ of $x$.
\end{enumerate}
\end{pro}

\pf
$(\Longrightarrow)$ Assume for simplicity that $\phi_t$ and $f_t$ are global
flows for $\xi$ and $x$. Since $\xi$ projects to $x$ under the groupoid
source projection $\alpha$, it follows that $\alpha\circ\phi_t =
f_t\circ\alpha$. Similarly $\beta\circ\phi_t = f_t\circ\beta$.

Denote $\{(h,g)\in G\times G \ST \alpha(h) = \beta(g)\}$ by $G*G$. Define a
vector field $\xi*\xi$ on $G*G$ by $\xi*\xi(h,g) = (\xi(h),\xi(g))$; since
$T(\alpha)(\xi(h)) = x(\alpha(h)) = x(\beta(g)) = T(\beta)(\xi(g))$,
$\xi*\xi$ is tangent to $G*G$. Evidently $\xi*\xi$ has flow $\psi_t(h,g) =
(\phi_t(h),\phi_t(g))$. Denoting the groupoid composition by $\kappa\colon
G*G\to G$, we know that $\xi*\xi$ projects to $\xi$ under $\kappa$. It follows
that $\phi_t(h)\phi_t(g) = \phi_t(hg)$ for all $(h,g)\in G*G$, and so $\phi_t$
is a Lie groupoid automorphism.

The converse is established by retracing these steps in the reverse order.
\qed

Call $F\colon G\to\R$ a {\em multiplicative function} if it is a groupoid
morphism from $G$ into the abelian group $\R$.

\begin{cor}
Let $(\xi,x)$ be a multiplicative vector field on $G$, and let
$F\colon G\to\R$ be a multiplicative function. Then $\xi(F)$ is also a
multiplicative function.
\end{cor}

It follows from Proposition~\ref{pro:mvf}, or from \cite[4.3]{Mackenzie:1992}, that if
$(\xi,x)$ and $(\eta,y)$ are multiplicative vector fields, then
$([\xi,\eta],\,[x,y])$ is also.

\begin{pro}
Let $(\xi,x)$ be a \mvf\ and take $X\in\Ga AG$ with corresponding right-
and left-invariant vector fields $\Vec{X}$ and $\ceV{X}$. Then $[\xi,\Vec{X}]$
is right-invariant and $[\xi,\ceV{X}]$ is left-invariant.
\end{pro}

\pf
Because $\xi$ is multiplicative we have, as in \ref{pro:mvf}, that
$\xi*\xi\stackrel{\kappa}{\sim}\xi$. Further, if $Z\in{\cal X}(G)$ is
$\alpha$-vertical, it is right-invariant iff $Z*0\stackrel{\kappa}{\sim}Z$
(see \cite[\S4]{Mackenzie:1992}). So we have
$[\xi*\xi,\Vec{X}*0]\stackrel{\kappa}{\sim}[\xi,\Vec{X}]$, whence
$[\xi,\Vec{X}]*0\stackrel{\kappa}{\sim}[\xi,\Vec{X}]$, and so $[\xi,\Vec{X}]$
is right-invariant.

A similar proof applies in the left-invariant case.
\qed

For a \mvf\ $(\xi,x)$ there is now a map
$$
D_\xi\colon\Ga AG\to \Ga AG,\qquad \Vec{D_\xi(X)} = [\xi,\Vec{X}].
$$
The proof of the following result is straightforward.

\begin{pro}
\begin{enumerate}
\item The map $D_\xi\colon\Ga AG\to \Ga AG$ is a derivation of the bracket
structure.

\item For $X\in\Ga AG$ and any vector field $\Bar{X}$ on $G$ with
$\Bar{X}\vert_M = X$, we have $D_\xi(X) = [\xi,\Bar{X}]\vert_M$.

\item If $(\eta,y)$ is a second \mvf, then $D_{[\xi,\eta]} = [D_\xi,D_\eta]$.
\end{enumerate}
\end{pro}

Given a multiplicative vector field $(\xi,x)$, applying the Lie functor
produces a section
\begin{equation}
\matrix{&ATG&\lrah&TM&\cr
        &&&&\cr
  A(\xi)&\Big\uparrow&&\Big\uparrow&x\cr
        &&&&\cr
        &AG&\lrah&M.&\cr
        &&q_G&&\cr}
\end{equation}

Define $\Tilde{\xi} = (j_G)^{-1}\circ A(\xi)$; it is clear that
$(\Tilde{\xi},x)$ is a linear vector field on $AG$. If $G = M\times M$ and
$\xi = x\times x$ where $x\in{\cal X}(M)$, then $\Tilde{\xi}$ is the
{\em complete lift} $\Tilde{x}$ of $x$ to $TM$ in the sense of \cite{YanoI}.
The linear vector field $(\Tilde{x},x)$ corresponds to the Lie derivative
$L_x\colon{\cal X}(M)\to{\cal X}(M),\ y\mapsto [x,y].$ From Proposition
\ref{prop:HigginsM:1990a} we therefore have
\begin{equation}                             \label{eq:brack}
\tau(x(m),[x,y](m)) = \Tilde{y}(x(m)) - T(x)(y(m)).
\end{equation}

If $\phi_t$ is a (local) flow for $\xi$, then $A(\phi_t)$ is a (local) flow
for $\Tilde{\xi}$; this follows from the classical result for complete
lifts to tangent bundles. Since $a_{TG}\circ j_G = J_M\circ T(a_G)$, the
vector field $\Tilde{\xi}$ on $AG$ projects under $a_G$ to $\Tilde{x}$ on
$TM$.

Since $(\Tilde{\xi},x)$ is linear, it induces a \cdo\
$D_{\Tilde{\xi}}\colon\Ga AG\to \Ga AG$ as in (\ref{eq:D}).

\begin{thm}                                    \label{thm:D}
For any \mvf\ $(\xi,x)$ on $G$, $D_{\Tilde{\xi}} = D_\xi.$
\end{thm}

\pf
Applying (\ref{eq:brack}) to $\xi$ and $\Vec{X}$ for $X\in\Ga AG$, at
$1_m\in G$, and recalling that $\xi(1_m) = T(1)(x(m))$, we have
$$
\tau(X(m),[\xi,\Vec{X}](1_m)) = T(\Vec{X})(T(1)(x(m)) - \Tilde{\xi}(X(m));
$$
note that $\tau$ and the tildes here refer to the double vector bundle $T^2G$.
On the other hand, in $TAG$ we have, from Proposition \ref{prop:HigginsM:1990a},
$$
\tau(X(m),D_{\Tilde{\xi}}(X)(m)) = T(X)(x(m)) - \Tilde{\xi}(X(m)).
$$
Note that $\Tilde{\xi}$ as a vector field on $AG$ is the restriction of
$\Tilde{\xi}$ as a vector field on $TG$.

>From \cite[\S7]{MackenzieX:1994} we know that
$$
\Tilde{\Vec{X}} = \Vec{j_G\circ T(X)},
$$
where the tilde on the left refers to $T^2G$ and the arrow on the right
refers to the tangent groupoid $TG$. So
$T(\Vec{X}) = J_G(\Vec{j_G\circ T(X)})$. Since $j_G\colon TAG\to ATG$ is
a restriction of $J_G\colon T^2G\to T^2G$, it follows that
$T(\Vec{X})(T(1)(x(m))) = T(X)(x(m)).$ Regarding $ATG\subseteq T^2G$ and
noting that the $\tau$ for $TAG$ is then the restriction of the $\tau$
for $T^2G$, we have $D_{\Tilde{\xi}}(X)(m) = [\xi,\Vec{X}](1_m).$
\qed

For any \mvf s $(\xi,x)$ and $(\eta,y)$ on $G$ and $f\in \cinf(M),\
\phi\in\Ga A^*G$, we now have
\begin{equation}                  \label{eq:Tildes}
\Tilde{\xi+\eta} = \Tilde{\xi} + \Tilde{\eta},\qquad
\Tilde{\xi}(f\circ q) = x(f)\circ q,\qquad
\Tilde{\xi}(\ell_{\phi}) = \ell_{D^{(*)}_\xi(\phi)}.
\end{equation}

For $F\colon G\to\R$ any multiplicative function, denote $A(F)\colon AG\to\R$
by $\Tilde{F}$.

\begin{pro}                             \label{pro:Tildes}
For $(\xi,x),(\eta,y)$ multiplicative vector fields on $G$, $X,Y\in\Ga AG$,
and $F$ a multiplicative function on $G$,
$$
{[\Tilde{\xi},\Tilde{\eta}]} = \Tilde{[\xi,\eta]},\qquad
{[\Tilde{\xi},X^\upa]} = D_\xi(X)^\upa,\qquad
{[X^\upa,Y^\upa]} = 0,\qquad
\Tilde{\xi}(\Tilde{F}) = \Tilde{\xi(F)}.
$$
\end{pro}

\pf
The third equation is known from (\ref{eq:cvf}), and the last follows from
Proposition~\ref{pro:mvf}.

For the first and second equations, it suffices to verify equality on
functions of the forms $\ell_\phi,\ \phi\in\Ga A^*G,$ and $f\circ q,\
f\in \cinf(M)$. For the second equation and $\phi\in\Ga A^*G$ we have
\begin{eqnarray*}
[\Tilde{\xi},X^\upa](\ell_\phi)
& = & \Tilde{\xi}(\langle\phi,X\rangle\circ q) -
X^\upa(\ell_{D^{(*)}_\xi(\phi)}) \\
& = & x(\langle\phi,X\rangle)\circ q -
\langle D^{(*)}_\xi(\phi),X\rangle\circ q \\
& = & \langle\phi,D_\xi(X)\rangle\circ q \\
& = & D_\xi(X)^\upa(\ell_\phi),
\end{eqnarray*}
whilst for $f\in \cinf(M)$,
$$
[\Tilde{\xi},X^\upa](f\circ q) = 0 = D_\xi(X)^\upa(f\circ q).
$$
The first equation is proved in a similar way.
\qed

Taking $G$ to be the pair groupoid $M\times M$ we of course recover the
classical description of vector fields on a tangent bundle \cite{YanoI}.

Multiplicative vector fields are rather special and we now briefly consider
two more general types of vector field on Lie groupoids. In what follows we
omit most proofs, which are similar to those given above.

\begin{defi}                               \label{defi:svf}
A {\em star vector field} on $G$ is a pair of vector fields $(\xi,x)$
where $\xi\in{\cal X}(G),\ x\in{\cal X}(M)$, such that
$T(\alpha)\circ\xi = x\circ\alpha$ and $\xi\circ 1 = T(1)\circ x$.
\end{defi}

The terminology comes from \cite{HigginsM:1990a}, where a {\em star map}
of groupoids is a map which preserves the $\alpha$-fibration and the
identities.

\begin{lem}                            \label{lem:svf}
Given any vector field $x$ on $M$, there is a star vector field
$(\xi,x)$ on $G$.
\end{lem}

\pf
Define a vector field $\eta$ on $G$ by setting $\eta(1_m) = T(1)(x(m))$ for
$m\in M$, and extending over $G$. Then $\mu\colon G\to TM$ defined by
$\mu(g) = T(\alpha)(\eta(g)) - x(\alpha g)$ is a section of the pullback
bundle $\alpha^!TM$. Since $\alpha$ is a surjective submersion, there is a
vector field $\zeta$ on $G$ with $T(\alpha)\circ\zeta = \mu$; we can also
require that $\zeta$ vanish on all $1_m\in G$. Now $\xi = \eta - \zeta$ is a
star vector field over $x$.
\qed

\begin{pro}                       \label{pro:svf}
Let $\xi$ be a vector field on a Lie groupoid $G$ and $x$ a vector field on
its base $M$. Then the following are equivalent:
\begin{enumerate}
\item $(\xi,x)$ is a star vector field on $G$;

\item The flows $\phi_t$ of $\xi$ are (local) star maps over the flows
$f_t$ of $x$.
\end{enumerate}
\end{pro}

It follows as before that if $(\xi,x)$ and $(\eta,y)$ are star vector fields,
then $([\xi,\eta],\,[x,y])$ is also. If $(\xi,x)$ is a star vector field and
$X\in \Ga AG$, then it is clear that $[\xi,\Vec{X}]$ is $\alpha$-vertical. We
now define $D_\xi\colon\Ga AG\to \Ga AG$ by
$$
D_\xi(X) = [\xi,\Vec{X}]\circ 1.
$$
Again, $D_\xi(X) = [\xi,\Bar{X}]\circ 1$ for any vector field $\Bar{X}$ on $G$
such that $\Bar{X}\circ 1 = X$.

Given a star vector field $(\xi,x)$, we can still apply the Lie functor and
obtain a linear vector field $\Tilde{\xi} = (j_G)^{-1}\circ A(\xi)$ on $AG$.
If $\phi_t$ is a (local) flow for $\xi$, then $A(\phi_t)$ is still a (local)
flow for $\Tilde{\xi}$. The proof of the following result follows as for
Theorem \ref{thm:D}.

\begin{thm}                                    \label{thm:svfD}
For any star vector field $(\xi,x)$ on $G$, $D_{\Tilde{\xi}} = D_\xi.$
\end{thm}

Equations (\ref{eq:Tildes}) continue to hold for star vector fields $(\xi,x)$
and $(\eta,y)$ and $f\in \cinf(M),\ \phi\in\Ga A^*G$.

\begin{pro}                             \label{pro:sTildes}
For $(\xi,x),(\eta,y)$ star vector fields on $G$, and $X,Y\in\Ga AG$,
$$
{[\Tilde{\xi},\Tilde{\eta}]} = \Tilde{[\xi,\eta]},\qquad
{[\Tilde{\xi},X^\upa]} = D_\xi(X)^\upa,\qquad
{[X^\upa,Y^\upa]} = 0.
$$
\end{pro}
In view of the following result, these equations determine the bracket
structure for all vector fields on $AG$.

\begin{pro}                          \label{pro:generation}
The vector fields of the form $\Tilde{\xi}$, where $(\xi,x)$ is a star vector
field, together with those of the form $X^\upa$, where $X\in\Ga AG$, generate
${\cal X}(AG)$.
\end{pro}

\pf
Take $\Xi\in TAG$ with $T(q)(\Xi) = x(m)$ and $p_{AG}(\Xi) = X$. Extend
$x(m)\in TM$ to a vector field $x$ on $M$. By Lemma \ref{lem:svf}, there is a
star vector field $(\xi,x)$ on $G$. We now have $T(q)(\Tilde{\xi}(X)) = x(m)$
and so, by the first sequence in (\ref{eq:pAcore}), we have
$$
\Xi = \Tilde{\xi}(X) + Y^\upa(X)
$$
for some $Y\in\Ga AG$.
\qed

Lastly in this section, we briefly consider the notion of affine vector
field on a Lie groupoid.

\begin{defi}
A vector field $\xi$ on $G$ is {\em affine} if for all $g,h\in G$ with
$\alpha(g) = \beta(h) = m$ we have
$$
\xi(gh) = T(L_{{\cal X}})(\xi(h)) + T(R_{{\cal Y}})(\xi(g)) -
T(L_{{\cal X}})T(R_{{\cal Y}})(\xi(1_m)),
$$
where ${\cal X}, {\cal Y}$ are bisections with ${\cal X}(\alpha g) = g,\
{\cal Y}(\alpha h) = h$.
\end{defi}

A multiplicative vector field is affine (see Proposition \ref{pro:sol}) and
for any $X\in\Ga AG$, both $\Vec{X}$ and $\ceV{X}$ are affine vector fields.
It is clear that the sum and scalar multiples of affine vector fields are
affine. Further, affine vector fields are closed under the bracket of
vector fields.

Any affine vector field on a Lie group is a sum of a multiplicative vector
field and a right-invariant vector field (see \cite[4.11]{Weinstein:1990} for
the case of bivector fields). Since
multiplicative and right-invariant vector fields on a groupoid are always both
$\alpha$- and $\beta$-projectable, the following result is the best that can
be expected in the general case. It is easy to construct affine vector fields
on a pair groupoid $M\times M$ which are not projectable.

\begin{pro}
Let $\xi$ be an affine vector field on the Lie groupoid $G$, and suppose
that $\xi$ is both $\alpha$- and $\beta$-projectable. Then $\xi$ is the
sum of a multiplicative and a right-invariant (or left-invariant) vector
field.
\end{pro}

\pf
Define $X\in\Ga AG$ by $X(m) = \xi(1_m) - T(1)T(\alpha)(\xi(1_m))$, where
$m\in M$, and write $\eta = \xi - \Vec{X}$. We prove that $\eta$ is
multiplicative. Let $x,y$ be the vector fields on $M$ for which
$\xi\stackrel{\alpha}{\sim}x,\ \xi\stackrel{\beta}{\sim} y$. Then $\Vec{X}
\stackrel{\beta}{\sim} y-x$ and so $\eta\stackrel{\alpha}{\sim}x$ and
$\eta\stackrel{\beta}{\sim} x$. Clearly $\eta(1_m) = T(1)(x(m))$ for
$m\in M$.

Since $\xi$ and $\Vec{X}$ are affine, it follows that $\eta $ is affine.
In this case, the affine condition implies that $\eta (gh) =
\eta (g) \sol \eta (h)$. Thus $\eta $ is multiplicative.
\qed

\section{Vector fields on Lie algebroid duals}        \label{sect:vfolad}

Now consider the Lie algebroid dual $q_*\colon A^*G\to M$ for a Lie groupoid
$G$ on base $M$. The tangent double vector bundle $(T(A^*G);A^*G,TM;M)$ has
core $A^*G$ and so for each $\phi\in\Ga A^*G$ there is a core vector field
$\phi^\upa\in{\cal X}(A^*G)$ such that
$$
\phi^\upa(\ell_X) = \langle\phi,X\rangle\circ q_*,\qquad\qquad
\phi^\upa(f\circ q_*) = 0,
$$
for all $X\in\Ga AG,\ f\in \cinf(M)$. Thus $[\phi^\upa,\psi^\upa] = 0$
for all $\phi,\psi\in\Ga A^*G$.

Let $(\xi,x)$ be a star vector field on $G$. Then the covariant differential
operators $D_\xi\in\Ga\CDO(AG)$ and $D^{(*)}_\xi\in\Ga\CDO(A^*G)$ induce a
linear vector field $(H_\xi,x)$ on $A^*G$. From (\ref{eq:D*}) and
the fact that $H_\xi\stackrel{q_*}{\sim} x$ we have
$$
H_\xi(\ell_X) = \ell_{D_\xi(X)},\qquad\qquad
H_\xi(f\circ q_*) = x(f)\circ q_*,
$$
for $X\in\Ga AG,\ f\in \cinf(M)$. The proof of the following result follows as
for Proposition \ref{pro:Tildes}.

\begin{pro}                               \label{pro:4.1}
For $(\xi,x),(\eta,y)$ star vector fields on $G$ and $\phi,\psi\in\Ga A^*G$,
$$
{[H_\xi,H_\eta]} = H_{[\xi,\eta]},\qquad
{[H_\xi,\phi^\upa]} = D_\xi^{(*)}(\phi)^\upa,\qquad
{[\phi^\upa,\psi^\upa]} = 0.
$$
\end{pro}

The next result shows in particular that for $(\xi,x)$ a \mvf\ on $G$,
the vector field $(H_\xi,x)$ is a Poisson vector field (called a Poisson
infinitesimal automorphism in \cite{LibermannM}) with respect to the
dual Poisson structure \cite{Courant:1990}, \cite{DazordS:1988} on $A^*G$.

\begin{pro}                          \label{pro:lpvf}
Let $(H,x)$ be a linear vector field on a Lie algebroid dual $A^*$, with
corresponding covariant differential operator $D\colon\Ga A\to\Ga A$. Then
$H$ is a Poisson vector field (that is, a derivation of the Poisson bracket
of functions on $A^*$) if and only if $D$ is a derivation of the bracket
in $\Ga A$.
\end{pro}

\pf Note first that if $D\in\Ga\CDO(A)$ is a derivation of the bracket in
$\Ga A$ then $a(D(X)) = [\Bar{a}(D),a(X)]$ for all $X\in\Ga A$, where
$\Bar{a}$ is the anchor on $\CDO(A)$; this follows by expanding the
expression $D([X,fY])$ in two ways.

Now the result is a straightforward calculation using the bracket relations
for the Poisson structure on $A^*$ \cite[equation (39)]{MackenzieX:1994}.
\qed

\begin{defi}
Let $A$ be a Lie algebroid on base $M$. Then a {\em morphic vector field}
on $A$ is a pair $(\Xi,x)$ where $\Xi\in{\cal X}(A),\ x\in{\cal X}(M)$,
such that $\Xi\colon A\to TA$ is a Lie algebroid morphism over
$x\colon M\to TM$.
\end{defi}

Here we take $TA$ as equipped with the tangent Lie algebroid structure on
base $TM$ described in \cite[\S5]{MackenzieX:1994}. Using
\cite[7.1]{MackenzieX:1994} it follows from the definition of $\Tilde{\xi}$
that if $(\xi,x)$ is a \mvf\ on a Lie groupoid $G$, then $(\Tilde{\xi},x)$
is a morphic vector field on $AG$.

\begin{thm}                          \label{thm:movf}
Let $A$ be a Lie algebroid on base $M$, and let $(\Xi,x)$ be a linear
vector field on $A$. The following are equivalent:
\begin{enumerate}
\item $(\Xi,x)$ is a morphic vector field;
\item The dual vector field $(\Xi_*,x)$ on $A^*$ is Poisson;
\item $D_\Xi\colon\Ga A\to\Ga A$ is a derivation;
\item The flows of $\Xi$ are (local) Lie algebroid automorphisms over the
flows of $x$.
\end{enumerate}
\end{thm}

\pf
The equivalence of (ii) and (iii) is Proposition \ref{pro:lpvf} above.
The equivalence of (ii) and (iv) follows from the facts that a vector
field is Poisson iff its flows are (local) Poisson automorphisms
\cite[III\S10]{LibermannM}, and that a vector bundle isomorphism of Lie
algebroids is a Lie algebroid isomorphism iff its dual is a Poisson
isomorphism. To prove the equivalence of (i) and (ii) we first need some
results of independent interest.

\begin{pro} {\rm \cite[4.6]{Xu:1995}}         \label{pro:subalg}
Let $B\lon S$ be a vector subbundle of a Lie algebroid $A\lon M$. Then $B$
is a Lie subalgebroid of $A$ if and only if $B^\perp$ is coisotropic in $A^*$
with the dual Poisson structure.
\end{pro}

\begin{pro}                                  \label{pro:phi}
Let $A$ be any Lie algebroid, and consider $\phi\in\Ga A^*$. Then $d\phi = 0$
if and only if $\im(\phi)$ is a coisotropic submanifold of $A^*$ with the dual
Poisson structure.
\end{pro}

\pf
For any $X\in \Ga (A)$, let $f_{X}$ be the function $\langle\phi,X\rangle$
on the base manifold $M$. It is clear that
$\ell_{X} -q^{*}f_{X}$ vanishes on the graph $\im \phi$. In fact, the space
of functions vanishing on $\im \phi$ is spanned by all such functions for
$X\in \Ga (A)$. Therefore, to show that $\im \phi$ is coisotropic, it
suffices to prove that $\{\ell_{X}-q^{*}f_{X}, \ \ell_{Y}-q^{*}f_{Y}\}$
vanishes on $\im \phi$ for any  $X, Y\in \Ga (A)$.

Now,
\be
\{\ell_{X}-q^{*}f_{X}, \ \ell_{Y}-q^{*}f_{Y}\}
&=&\{\ell_{X} , \ell_{Y}\}-\{q^{*}f_{X}, \ell_{Y}\}-\{\ell_{X}, q^{*}f_{Y}\}
   +\{q^{*}f_{X}, q^{*}f_{Y}\}\\
&=&\ell_{[X, Y]}+q^{*}(a(Y)f_{X})-q^{*}(a(X)f_{Y}),
\ee
where the last step follows directly from the definition of the Poisson
structure on $A^*$. By evaluating at $\phi (m)$, one obtains that
\be
\{\ell_{X}-q^{*}f_{X}, \ \ell_{Y}-q^{*}f_{Y}\} (\phi (m ))
&=&\langle [X, Y], \phi \rangle (m) + (a(Y)\langle X, \phi \rangle )(m)
   - (a(X)\langle Y, \phi \rangle )(m)\\
&=&-(d\phi )(X, Y)(m).
\ee

This leads to our conclusion: $\im \phi $ is coisotropic iff $\phi $ is
closed.
\qed

\begin{rmk}\rm
When $A$ is a tangent bundle Lie algebroid $TM$, $A^*$ is the cotangent bundle
$T^*M$ with the canonical symplectic structure. In this case Proposition
\ref{pro:phi} is just the well-known fact,
that the graph of a one-form is Lagrangian iff the form is closed. On the
other hand, if $A$ is a Lie algebra $\frakg$, a one-form  $\phi$
is just a point in $\frakg^*$. It is coisotropic iff the Poisson
tensor vanishes at the point, which is exactly what $d\phi =0$ means.
\end{rmk}

\begin{pro}                                   \label{pro:pvf}
Let $X$ be a vector field on a Poisson manifold $P$. Then $X$ is Poisson
if and only if $\im(X)$ is a coisotropic submanifold of $TP$ with the
tangent Poisson structure.
\end{pro}

\pf Take $A = T^*P$ in the preceding Proposition. Then for $X\in\Ga TP =
\Ga A^*$ we have $dX = [\pi,X] = L_X(\pi)$, where $\pi$ is the Poisson
tensor. Thus $dX = 0$ if and only if $X$ is a Poisson vector field. The
result now follows from the fact that the dual Poisson structure on $TP$
from the Lie algebroid structure on $T^*P$ is the tangent Poisson structure.
\qed

We can now complete the proof of Theorem \ref{thm:movf}. Firstly, since $\Xi$
and $x$ are embeddings, $\Xi\colon A\lon TA$ is a morphism over $x\colon M\lon
TM$ if and only if $\im(\Xi)\lon\im(x)$ is a Lie subalgebroid of $TA\lon TM$.

Next recall from \cite[\S5]{MackenzieX:1994} the dual $T^\sol(A)\lon TM$ of
the vector bundle $TA\lon TM$. The tangent pairing (\ref{eq:pairing}) induces
an isomorphism $I\colon T(A^*)\lon T^\sol(A)$, and this is also an isomorphism
of the Poisson structures on $T(A^*)$ and $T^\sol(A)$
\cite[5.6]{MackenzieX:1994}. It now follows easily from Proposition \ref{xiH}
that $I$ maps $\im(\Xi_*)$ isomorphically onto $\im(\Xi)^\perp$. Using
Propositions \ref{pro:subalg} and \ref{pro:pvf}, the proof of Theorem
\ref{thm:movf} is complete.
\qed

Theorem \ref{thm:movf} is slightly surprising, because it is certainly not
true that a vector field on a Poisson manifold is Poisson iff it is a
Poisson morphism into the tangent Poisson structure.

\begin{thm}                                 \label{thm:intvf}
Let $G$ be an $\alpha$-simply connected Lie groupoid on base $M$. Then
if $(\Xi,x)$ is a morphic vector field on $AG$, there is a unique \mvf\
$(\xi,x)$ on $G$ such that $\Tilde{\xi} = \Xi$.
\end{thm}

\pf
This is an immediate consequence of the above and the integrability of
Lie algebroid morphisms \cite{MackenzieX2}.
\qed

\section{Forms on Lie groupoids}               \label{sect:flg}

We begin by considering any vector bundle $(A,q,M)$, and the double vector
bundle structure $(T^*(A);A,A^*;M)$ described in \cite[\S5]{MackenzieX:1994}.
Using the structure of $T^*A$, one can develop a concept of linear
differential form on $A$ corresponding to the concept of linear vector field
in \S\ref{sect:vfolg}. Here we will only present the few facts which we need.

\begin{defi}                             \label{df:l1f}
A {\em linear 1-form} on $A$ is a pair $(\Upsilon,\phi)$ where
$\Upsilon\in\fms{A}$ and $\phi\in\Ga A^*$ such that
\begin{equation}                       \label{diag:T*A}
\matrix{&&r&&\cr
        &T^*A&\lrah&A^*&\cr
        &&&&\cr
\Upsilon&\Big\uparrow&&\Big\uparrow&\phi\cr
        &&&&\cr
        &A&\lrah&M&\cr
        &&   &&\cr}
\end{equation}
is a morphism of vector bundles.
\end{defi}

Here the map $r$ is defined by $\langle r(\Phi),X\rangle =
\langle\Phi,X^\upa(Y)\rangle$ for $\Phi\in T^*_YA$.

Note that the base map of a linear 1-form is not a 1-form on the base but
rather a section of the dual bundle. The sum of linear 1-forms is a linear
1-form and scalar multiples of linear 1-forms are linear 1-forms. The pairing
of a linear vector field with a linear 1-form is a fibrewise linear function.

Now consider a Lie groupoid $G$ on $M$. The cotangent bundle $T^*G$ has a
natural Lie groupoid structure on base $A^*G$ defined in \cite{CDW} and
\cite{Pradines:1988}. As in \cite[\S7]{MackenzieX:1994}, we take the source
$\tilalpha$ and target $\tilbeta$ to be given by
\begin{equation}                    \label{eq:T*G}
\Tilde{\alpha}(\omega)(X) = \omega(T(L_g)(X - T(1)(a(X)))),\qquad
\Tilde{\beta}(\omega)(Y) = \omega(T(R_g)(Y)),
\end{equation}
where $\omega\in T^*_gG,\ X\in A_{\alpha g}G$ and $Y\in A_{\beta g}G$. If
$\theta\in T^*_hG$ and $\Tilde{\alpha}(\theta) = \Tilde{\beta}(\omega)$
then $\alpha h = \beta g$ and we define $\theta\sol\omega\in T^*_{hg}G$ by
$$
(\theta\sol\omega)(Y\sol X) = \theta(Y) + \omega(X),
$$
where $Y\in T_hG,\ X\in T_gG$. The identity element
$\Tilde{1}_\phi \in T^*_{1_m}G$ corresponding to $\phi\in A^*_mG$ is defined
by $\Tilde{1}_\phi(T(1)(x) + X) = \phi(X)$ for $X\in A_mG,\ x\in T_m(M)$.

\begin{defi}                                 \label{df:m1f}
A {\em multiplicative 1-form} on $G$ is a pair $(\Phi,\phi)$ where
$\Phi\in\fms{G}$ and $\phi\in\Ga A^*G$ such that
\begin{equation}
\matrix{&T^*G&\lgpd&A^*G&\cr
        &&&&\cr
    \Phi&\Big\uparrow&&\Big\uparrow&\phi\cr
        &&&&\cr
        &G&\lgpd&M&\cr
        &&   &&\cr}
\end{equation}
is a morphism of Lie groupoids.

A {\em star 1-form} on $G$ is similarly a pair $(\Phi,\phi)$ where
$\Phi\in\fms{G}$ and $\phi\in\Ga A^*G$ such that
$\tilalpha\circ\Phi = \phi\circ\alpha$ and $\Phi\circ 1 = \tilone\circ\phi$.
\end{defi}

Suppose $(\Phi,\phi)$ is a multiplicative 1-form and $(\xi,x)$ is a
multiplicative vector field on $G$. Writing
$F = \langle\Phi,\xi\rangle\colon G\to\R$, we have, for $h,g\in G$
compatible,
$$
F(hg) = \langle\Phi(h)\sol\Phi(g),\xi(h)\sol\xi(g)\rangle
$$
and it follows directly from the definition of the multiplication $\sol$
in $T^*G$ that this is equal to $\langle\Phi(h),\xi(h)\rangle +
\langle\Phi(g),\xi(g)\rangle.$ Thus $F$ is a multiplicative function.

Put differently, the standard pairing $T^*G\pback{G}TG\to\R$ is a groupoid
morphism with respect to the pullback groupoid structure with base
$A^*G\pback{M}TM$. Applying the Lie functor as in \cite[7.2]{MackenzieX:1994}, we obtain
a pairing $\llangle\ ,\ \rrangle\colon AT^*G\pback{AG}ATG\to\R$. This yields
an isomorphism of double vector bundles $i_G\colon AT^*G\to A^\sol TG$ defined
by
$$
\langle i_G({\goth X}),\xi\rangle^\sol = \llangle{\goth X},\xi\rrangle
$$
for ${\goth X}\in AT^*G,\ \xi\in ATG$. Here $A^\sol TG$ is the dual of $ATG$
over $AG$ and $\langle\ ,\ \rangle^\sol$ is the standard pairing of
$A^\sol TG$ and $ATG$. Still following \cite{MackenzieX:1994}, take the dual
of $j_G\colon TAG\lon ATG$ over $AG$ and define $j_G' = j_G^*\circ
i_G\colon AT^*G\to T^*AG$; this $j_G'$ is an isomorphism of double vector
bundles preserving $AG, A^*G$ and the cores. The following result is now
immediate.

\begin{pro}
For ${\goth X}\in AT^*G$ and $\xi\in TAG$,
$$
\llangle{\goth X},j_G(\xi)\rrangle = \langle j_G'({\goth X}),\xi\rangle.
$$
\end{pro}

Now consider a star 1-form $(\Phi,\phi)$ on $G$. Applying the Lie
functor gives a vector bundle morphism $A(\Phi)\colon AG\to AT^*G$. Define
$\Tilde{\Phi} = j_G'\circ A(\Phi)$; since $j_G'$ is an isomorphism of
double vector bundles, $(\Tilde{\Phi},\phi)$ is a linear 1-form on $AG$.

\begin{pro}                                \label{pro:s1f}
\begin{enumerate}
\item Let $(\Phi,\phi)$ be a star 1-form and $X\in\Ga AG$. Then
\begin{equation}                           \label{eq:tilphixupa}
\langle\Tilde{\Phi},X^\upa\rangle = \langle\phi,X\rangle\circ q.
\end{equation}

\item Let $(\Phi,\phi)$ be a star 1-form and $(\xi,x)$ a star vector field on
$G$. Then
\begin{equation}                            \label{eq:tilphitilxi}
\langle\Tilde{\Phi},\Tilde{\xi}\rangle = \Tilde{\langle\Phi,\xi\rangle},
\end{equation}
\nin where for any function $F\colon G\lon\R$, we denote
$A(F)\colon AG\lon\R$ by $\Tilde{F}$.
\end{enumerate}
\end{pro}

\pf
(i) follows from the fact that $r\circ\Tilde{\Phi} = \phi\circ q$, as in
(\ref{diag:T*A}). For (ii), since
$\llangle\ ,\ \rrangle = A(\langle\ ,\ \rangle)$, we have
$$
\langle\Tilde{\Phi},\Tilde{\xi}\rangle =
\langle j_G'\circ A(\Phi),j_G^{-1}\circ A(\xi)\rangle =
\llangle A(\Phi),A(\xi)\rrangle =
A(\langle\Phi,\xi\rangle) =
\Tilde{\langle\Phi,\xi\rangle}.
$$
\qed

These two equations describe the behaviour of the forms
$(\Tilde{\Phi},\phi)$ on ${\cal X}(AG)$. To complete the description of
the 1-forms on $AG$ we need to include the pullbacks of forms on the base
manifold.

Given $\omega\in\fms{M}$ there is the 1-form $q^*\omega\in\fms{AG}$.
Here it is immediate that
\begin{equation}                 \label{eq:pullbacks}
\langle q^*\omega,\Tilde{\xi}\rangle = \langle\omega,x\rangle\circ q,
\qquad\qquad
\langle q^*\omega,X^\upa\rangle = 0,
\end{equation}
for star vector fields $(\xi,x)$ on $G$ and $X\in\Ga AG$, since
$\Tilde{\xi}$ projects to $x$ under $q$ and $X^\upa$ projects to 0.

\begin{lem}                            \label{lem:sf}
Given any $\phi\in\Ga A^*G$, there is a star 1-form $(\Phi,\phi)$ on $G$.
\end{lem}

\pf
Define a 1-form $\Psi$ on $G$ by setting $\Psi(1_m) = \tilone_{\phi(m)}$ for
$m\in M$, and extending over $G$. Then $\mu\colon G\to A^*G$ defined by
$\mu(g) = \tilalpha(\Psi(g)) - \phi(\alpha g)$ is a section of the pullback
bundle $\alpha^!A^*G$. Since $\tilalpha\colon T^*G\lon A^*G$ is a fibrewise
surjection, there is a 1-form $\Upsilon$ on $G$ with $\tilalpha\circ\Upsilon
= \mu$; we can also require that $\Upsilon$ vanish on all $1_m\in G$. Now
$\Phi = \Psi - \Upsilon$ is a star 1-form over $\phi$.
\qed

\begin{pro}                          \label{pro:genforms}
The 1-forms $\Tilde{\Phi}$, where $(\Phi,\phi)$ is a star 1-form on $G$,
together with the pullbacks $q^*\omega$, where $\omega\in\fms{M}$, generate
$\fms{AG}$.
\end{pro}

\pf
Take $\Upsilon\in T^*_XAG$ with $r(\Upsilon) = \phi(m)$. Extend $\phi(m)$ to
a section $\phi$ of $\Ga A^*G$. By Lemma \ref{lem:sf}, there is a star 1-form
$(\Phi,\phi)$ on $G$. So at $X$ we have $r(\Tilde{\Phi}(X)) = \phi(m)$
and so,
$$
\Upsilon = \Tilde{\Phi}(X) + q^*\omega
$$
for some $\omega\in\fms{M}$.
\qed

Lastly in this section, we briefly indicate how to characterize the
$\Tilde{\Phi}$ corresponding to multiplicative 1--forms. For any Lie
algebroid $A$ on $M$, the dual $A^*$ has its Poisson structure and
$T^*A^*\to A^*$ has therefore a cotangent Lie algebroid structure (see
\cite{MackenzieX:1994} for detail and references). For any vector bundle
$A$ there is a canonical map $R\colon T^*A^*\to T^*A$ which is an isomorphism
of double vector bundles preserving $A$ and $A^*$ \cite{MackenzieX:1994}.
Using $R$ we transfer the Lie algebroid structure of $T^*A^*\to A^*$ to
$T^*A\to A^*$.

Now consider a multiplicative 1--form $(\Phi,\phi)$ on a Lie groupoid
$G\gpd M$. We have $\Tilde{\Phi} = j'_G\circ A(\Phi)$ and by
\cite[7.3]{MackenzieX:1994} we know $j'_G = R\circ s^{-1}$, where
$s\colon T^*A^*G\to AT^*G$ is the canonical isomorphism between the Lie
algebroid of the symplectic groupoid $T^*G$ and the cotangent Lie algebroid
of its base. Since $A(\Phi)$ and $s$ are Lie algebroid morphisms, it
follows that $\Tilde{\Phi}$ is a morphic 1--form on $AG$ in the sense of the
following definition.

\begin{defi}                            \label{df:mor1f}
Let $A$ be a Lie algebroid on $M$. A linear 1--form $(\Upsilon,\phi)$ on
$A$ is a {\em morphic 1--form} if $\Upsilon\colon A\to T^*A$ is a Lie
algebroid morphism over $\phi\colon M\to A^*$.
\end{defi}

The next result now follows just as in Theorem \ref{thm:intvf}.

\begin{thm}                                 \label{thm:int1f}
Let $G$ be an $\alpha$-simply connected Lie groupoid on base $M$. Then
if $(\Upsilon,\phi)$ is a morphic 1--form on $AG$, there is a unique
multiplicative 1--form $(\Phi,\phi)$ on $G$ such that
$\Tilde{\Phi} = \Upsilon$.
\end{thm}

\section{Forms on Lie algebroid duals}  \label{sect:folad}

Before we consider forms, notice that there is another way in which one may
extend the calculus of vertical and complete lifts to a Lie algebroid.
Consider any Lie algebroid $A$ on base $M$. Each $X\in\Ga A$ defines an
element $L_X$ of $\Ga\CDO(A)$ by $L_X(Y) = [X,Y]$ for which the corresponding
element of $\Ga\CDO(A^*)$ is the Lie derivative, also denoted $L_X$, defined
by
$$
L_X(\phi)(Y) = a(X)(\phi(Y)) - \phi([X,Y])
$$
for $\phi\in\Ga A^*$. There is a corresponding morphic vector field
$(\Tilde{X},a(X))$ on $A$. By definition, we now have
\begin{equation}
\Tilde{X}(\ell_\phi) = \ell_{L_X(\phi)},\qquad\qquad
\Tilde{X}(f\circ q) = a(X)(f)\circ q,
\end{equation}
for $\phi\in\Ga A^*,\ f\in C(M)$. Similarly, there is a linear Poisson
vector field $(H_X,a(X))$ on $A^*$ such that
\begin{equation}               \label{eq:HX}
H_X(\ell_Y) = \ell_{[X,Y]},\qquad\qquad
H_X(f\circ q_*) = a(X)(f)\circ q_*,
\end{equation}
for $Y\in\Ga A,\ f\in C(M)$. The following result is now immediate.

\begin{pro}            \label{other}
For $X,Y\in\Ga A$, and $\phi, \psi \in\Ga A^*$,
\begin{eqnarray*}
{[\Tilde{X},\Tilde{Y}]} = \Tilde{[X,Y]}, \qquad
{[\Tilde{X},Y^\upa]} = [X,Y]^\upa, \qquad
{[X^\upa,Y^\upa]} = 0, \\
{[H_X,H_Y]} = H_{[X,Y]}, \qquad
{[H_X,\phi^\upa]} = L_X(\phi)^\upa, \qquad
{[\phi^\upa,\psi^\upa]} = 0, \\
\end{eqnarray*}
\end{pro}

Although this formalism is valid for any Lie algebroid, not necessarily
integrable, the vector fields $\Tilde{X}$ and $X^\upa$ do not generally
generate the module of vector fields on $A$. For example, if $A$ is totally
intransitive and abelian (that is, a vector bundle) then the vector fields
$\Tilde{X}$ are all zero.

We now describe forms on a Lie algebroid dual in terms of forms of two
specific types.
Whereas linear vector fields on the dual of a vector bundle $A$ correspond
to linear vector fields on $A$, there is no such result for 1--forms, and
the description we give here is not related to that for 1--forms on $A$ in
the same way that our descriptions of vector fields on $AG$ and $A^*G$ are
related. In fact, the following
description applies to any vector bundle.

Firstly, given any $\omega\in\fms{M}$, there is the corresponding 1-form
$q_*^*\omega\in\fms{A^*}$. Evidently, for $\phi\in\Ga A^*$ and
$X\in\Ga A$,
$$
\langle q_*^*\omega,\phi^\upa\rangle = 0,\qquad\qquad
\langle q_*^*\omega,H_X\rangle = \langle\omega,a(X)\rangle\circ q_*.
$$
If $A = AG$ is integrable and $(\xi,x)$ is a star vector field on $G$,
we also have
$$
\langle q_*^*\omega,H_\xi\rangle = \langle\omega,x\rangle\circ q_*.
$$

Secondly, for each $X\in\Ga A$, there is the form $\delta\ell_X\in\fms{A^*}$.
>From (\ref{eq:HX}) it follows that, for $X,Y\in\Ga A$ and
$\phi\in\Ga A^*$,
$$
\langle\delta\ell_X,H_Y\rangle = \ell_{[Y,X]},\qquad\qquad
\langle\delta\ell_X,\phi^\upa\rangle = \langle\phi,X\rangle\circ q_*.
$$
If $A = AG$ is integrable and $(\xi,x)$ is a star vector field on $G$,
we also have
$$
\langle\delta\ell_X,H_\xi\rangle = \ell_{D_\xi(X)}.
$$

By a result similar to Proposition \ref{pro:genforms}, forms of these two
types generate $\fms{A^*}$.

Now $A^*$ has the Poisson structure dual to its Lie algebroid structure.
By \cite[5.7]{MackenzieX:1994}, the Poisson anchor
$\pi^\#\colon T^*(A^*)\to T(A^*)$
is a morphism of double vector bundles, as in Figure~\ref{fig:pa},
\begin{figure}[htb]
\begin{picture}(350,200)(-100,0)
\put(-70,155){$\matrix{&&       &\cr
                      &T^*(A^*)&\lrah &A\cr
                      &&&\cr
                      &\Bigg\downarrow&&\Bigg\downarrow\cr
                      &&&\cr
                      &A^*&\lrah&M\cr}$}
\put(-30,170){\vector(1,-1){60}}            \put(10,150){$\pi^\#$}
\put(70,170){\vector(1,-1){60}}            \put(100,150){$a$}
\put(10,70){$\matrix{&&      &\cr
                     &T(A^*)&\lrah &TM\cr
                     &&&\cr
                     &\Bigg\downarrow&&\Bigg\downarrow\cr
                     &&&\cr
                     &A^*&\lrah&M\cr}$}
\end{picture}
\caption{\ \label{fig:pa}}
\end{figure}
and the induced map of the cores is $-a^*\colon T^*M\to A^*$.

For ${\goth F}\in\fms{A^*}$, denote the vector field on $A^*$ corresponding
under $\pi^\#$ by ${\goth F}^\#$. Then we have, in particular, that for
$\omega\in\fms{M}$,
$$
(q_*^*\omega)^\# = -a^*(\omega)^\upa;
$$
that is, the Hamiltonian vector field corresponding to the pullback form
$q_*^*\omega$ is the negative of the vertical lift corresponding to
$a^*(\omega)\in\Ga A^*$.

Further, it is easy to check that, for $X\in\Ga A$,
$$
(\delta\ell_X)^\# = H_X.
$$
Thus $H_X$ is the Hamiltonian vector field with energy $\ell_X$.

\section{Fields and forms on Poisson groupoids}
\label{sect:fafopg}

In \cite{Weinstein:1988}, Weinstein showed that the Lie algebroid dual of a
Poisson groupoid $G$ has itself a Lie algebroid structure; he obtained this
by means of a general linearization result along an arbitrary coisotropic
submanifold of any Poisson manifold. In \cite{MackenzieX:1994} we obtained
this structure by applying the Lie functor to the Poisson tensor
$\pi^\#\colon T^*G\lon TG$ and showing that, after certain canonical
isomorphisms, this defines a Poisson structure on $AG$ whose dual is the
desired Lie algebroid structure on $A^*G$. Here we provide a particularly
concrete and explicit description of this structure. Firstly, in Theorem
\ref{thm:last} we show how the bracket of 1-forms on the Poisson manifold
$AG$ is determined in terms of the bracket of 1-forms of $G$; this uses the
process $\Phi\mapsto\Tilde{\Phi}$ of \S\ref{sect:flg} and a family of
covariant differential operators $D_\Phi$ in $T^*P$, induced by the star
1-forms $\Phi$ on $G$. Secondly, equation (\ref{eq:lbabrack}) is an explicit
formula for the bracket on $\Ga A^*G$ in terms of the anchor of $A^*G$ and
the Poisson tensor $T^*G\to TG$ (note that, by Lemma \ref{lem:sf}, every
section of $A^*G$ is the base of some star 1-form on $G$). If $G$ is a pair
(or coarse) Poisson groupoid $\Bar{P}\times P$, then (\ref{eq:lbabrack})
reduces to the usual formula for the bracket of 1-forms on a Poisson
manifold.

Throughout this section, $G$ denotes a Poisson groupoid on base $P$. We
follow the notation and conventions of \cite{MackenzieX:1994}.

Consider a star 1-form $(\Phi,\phi)$ on $G$. Analagously to the construction
in \S\ref{sect:vfolg}, we define an operator $D_\Phi$ on $T^*P$. Note first
that for any $\omega\in\fms{P}$, the pullback $\beta^*\omega\in\fms{G}$
projects under $\tilalpha\colon T^*G\lon A^*G$ to the zero section of
$A^*G$. Since $\Phi$ projects to $\phi$ under $\tilalpha$ by assumption, and
$\tilalpha$ is a Lie algebroid morphism, it follows that
$[\Phi,\beta^*\omega]$ also projects to zero; that is,
\begin{equation}                            \label{eq:Phi}
\langle[\Phi,\beta^*\omega],T(L_g)(X-T(1)(aX))\rangle = 0
\end{equation}
for all $g\in G$ and $X\in A_{\alpha g}G$. Now every $\beta$-vertical vector
is the left-translate of some $Y\in T_{1_m}(\beta^{-1}(m))$ and every such
$Y$ is equal to $T(i)(X) = -X + T(1)(aX)$ for some $X\in A_mG$. It therefore
follows from (\ref{eq:Phi}) that $[\Phi,\beta^*\omega]$ annuls all
$\beta$-vertical vectors. We can therefore define $D_\Phi(\omega)\in\fms{P}$
by
\begin{equation}                            \label{eq:defn}
\langle D_\Phi(\omega),T(\beta)(Y)\rangle =
         \langle[\Phi,\beta^*\omega],Y\rangle
\end{equation}
where $Y\in T_{1_m}G$.

\begin{pro}
\begin{enumerate}
\item For $(\Phi,\phi)$ a star 1-form, $D_\Phi\colon\fms{P}\lon\fms{P}$
is a covariant differential operator over $a_*(\phi)$; that is, it is
$\R$-linear and $D_\Phi(f\omega) = fD_\Phi(\omega) + a_*(\phi)(f)\omega$
for all $f\in\cinf(P),\ \omega\in\fms{P}$;

\item for any $Y\in T_{1_m}G$ and any $\theta\in\fms{G}$ such that
$\theta(1_m) = \beta^*\omega(1_m)$ for all $m\in P$,
$$
\langle D_\Phi(\omega),T(\beta)(Y)\rangle =
         \langle[\Phi,\theta],Y\rangle;
$$

\item if $(\Psi,\psi)$ is another star 1-form, then
$D_{[\Phi,\Psi]} = [D_\Phi,D_\Psi].$
\end{enumerate}
\end{pro}

\pf
The proof of (i) is straightforward, and (iii) follows directly from (ii).
To prove (ii), it is sufficent to show that if $\theta\in\fms{G}$ vanishes
on the identity elements of $G$, then $[\Phi,\theta]$ does also. But
$\tilone\colon A^*G\lon T^*G$ is known to be a Lie algebroid morphism, so
$0\stackrel{1}{\sim}\theta$ and $\phi\stackrel{1}{\sim}\Phi$ imply that
$0 = [\phi,0]\stackrel{1}{\sim}[\Phi,\theta].$
\qed

\begin{thm}
Let $(\Phi,\phi)$ be a multiplicative 1-form on $G$. Then for all
$\omega\in\fms{P}$,
\begin{enumerate}
\item $a^*D_\Phi(\omega) = [\phi,a^*\omega]$,
\item $[\Phi,\beta^*\omega] = \beta^*(D_\Phi(\omega))$,
\item $D_\Phi$ is a derivation of the Poisson bracket on $\fms{P}$.
\end{enumerate}
\end{thm}

\pf
(i) For any $X\in \Ga (AG)$, using (\ref{eq:defn}),
$$
\langle a^{*}D_{\Phi}\omega , X\rangle  =
\langle D_{\Phi}\omega,aX\rangle = \langle D_{\Phi}\omega,T\beta X\rangle
= \langle [\Phi , \beta^{*} \omega ], X\rangle,
$$
where we consider $X$ as a section of $T_{P}^{\alpha }G $.

As a section of $A^*G,\ a^{*}\omega$ can be considered as a conormal one-form
along the base space $P$. Let $\widetilde{a^{*}\omega}$ be a one-form on $G$
extending $a^{*}\omega$. Now it suffices to show that
$\langle [\Phi , \beta^{*} \omega ], X\rangle =
\langle [\Phi ,  \widetilde{a^{*}\omega} ], X\rangle $, since the latter is, by
definition (see \cite{Weinstein:1988}),
$\langle [\phi , a^{*}\omega ], X\rangle $. Let
$\widetilde{\omega}=\beta^{*} \omega -\widetilde{a^{*}\omega}$. Then it
is clear that $\widetilde{\omega}=0$ when restricted to $T_{P}^{\alpha }G $.
Therefore, $\widetilde{\beta}(\widetilde{\omega }|_{P})=0$.
Since $\widetilde{\beta} : T^{*}G\lon A^{*}G$ is a Lie algebroid morphism
over $\beta : G\lon P$, we have
$\widetilde{\beta}([\Phi , \widetilde{\omega }]|_{P})=0$. This is to say
that $[\Phi , \beta^{*} \omega ] =
[\Phi ,  \widetilde{a^{*}\omega} ]$ when restricted to
$T_{P}^{\alpha }G$. This concludes the proof for (i).

(ii). We need a lemma before proving (ii). Let
$\Lambda \subset G\times G\times G$ be the graph of the groupoid
multiplication.

\begin{lem}
\label{lem:pullback}
Suppose that $\theta $ is a one-form on $G$. Then $\theta $ is the
pull back of a one-form on the base manifold $P$ iff
$(0, \theta , -\theta ) $ is conormal to $\Lambda $.
\end{lem}

\pf
Suppose that $(0, \theta , -\theta ) $ is conormal to $\Lambda $.
Let $\delta_{y}$ and $\delta_{z}$ be any tangent vectors
on $G$ such that $\beta_{*} \delta_{y}=\beta_{*} \delta_{z}$.
This condition, of course, implies  that $\beta (y)=\beta (z)$.
Let $x=zy^{-1}$. Then there exists a vector $\delta_{x}\in T_{x}G$
such that $(\delta_{x} , \delta_{y}, \delta_{z} )$ is tangent
to the graph $\Lambda $. It thus follows, by assumption, that
$\langle \theta , \delta_{y}\rangle =\langle \theta , \delta_{z}\rangle $, which
implies immediately that $\theta =\beta^{*}\psi $ for some
one-form $\psi $ on the base manifold.
The other direction is obvious.
\qed

Since $\Phi : G\lon T^{*}G $ is a groupoid morphism, the triple product
$\Phi \times \Phi \times \Phi$ maps $\Lambda$ into the graph of the groupoid
multiplication of $T^{*}G$. The latter is the space of all
$(\xi ,\eta , \zeta) \in T^{*}G \times T^{*}G \times T^{*}G$ such that
$(\xi ,\eta , - \zeta)$ is conormal to $\Lambda$.
In other words, this is equivalent to saying that
$(\Phi , \Phi, -\Phi )$,  as a one-form on $G\times G\times G$,
is conormal to $\Lambda$. On the other hand, we know
that $(0, \beta^{*}\omega , -\beta^{*} \omega )$ is conormal
to $\Lambda$ according to the lemma above. Since $\Lambda$ is a coisotropic
submanifold of $G\times G\times \Bar{G}$, we have that
$(0, [\Phi ,  \beta^{*} \omega ], -[\Phi ,  \beta^{*} \omega ])$
is conormal to $\Lambda$. It thus follows from  Lemma \ref{lem:pullback}
that $ [\Phi ,\beta^{*} \omega ] = \beta^{*}\psi $ for some one-form
$\psi$ on the base $P$. Since $\beta$ is a submersion, it is now clear that
$\psi = D_{\Phi} ( \omega )$.

(iii) Let $\omega, \ \theta $ be any one-forms on $P$,
and take $Y \in T_{1_{m}}G$. Using (ii),
\be
\langle D_{\Phi}[\omega , \theta ], T\beta Y\rangle
& = & \langle \beta^{*}D_{\Phi }[\omega , \theta ],Y\rangle \\
& = &\langle [\Phi , \beta^{*}[\omega , \theta ]], Y\rangle \\
& = &-\langle [\Phi , [\beta^{*}\omega , \beta^{*}\theta ]], Y\rangle
\ \ (\mbox{recall } \beta : G\lon P\ \mbox{ is anti-Poisson})\\
&=&-\langle [[ \Phi , \beta^{*}\omega ], \beta^{*}\theta ], Y\rangle
   -\langle [\beta^{*}\omega , [ \Phi , \beta^{*}\theta ]], Y\rangle \\
&=&-\langle [\beta^{*}  D_{\Phi} \omega , \beta^{*}\theta ], Y\rangle -
\langle [\beta^{*} \omega , \beta^{*}  D_{\Phi}\theta ], Y\rangle \\
&=&\langle  [ D_{\Phi} \omega , \theta ], T\beta Y\rangle +
   \langle [\omega , D_{\Phi}\theta ], T\beta Y\rangle . \\
\ee
It thus follows that $D_{\Phi}[\omega , \theta ]=[ D_{\Phi} \omega , \theta ]
+[\omega , D_{\Phi}\theta ]$.
\qed

Since $\pi^\#\colon T^*G\lon TG$ is a groupoid morphism over $a_*\colon
A^*G\lon TP$, it follows that if $(\Phi,\phi)$ is a star 1-form, then
$(\Phi^\#,a_*(\phi))$ is a star vector field, and if $(\Phi,\phi)$ is
multiplicative, then $(\Phi^\#,a_*(\phi))$ is also. Further, it is easily
seen that
$$
a_*^*(D_\Phi(\omega)) = D_{\Phi^\#}(a^*_*(\omega))
$$
for all $\omega\in\fms{P}$.

\begin{lem}
Let $(\Phi,\phi)$ be a star 1-form on $G$, and write $\xi = \Phi^\#$.
Then $\Tilde{\Phi}^\# = \Tilde{\xi}$ with respect to the Poisson structure on
$AG$ induced by the Lie algebroid structure on $A^*G$.
\end{lem}

\pf
Apply $A$ to $\xi = \pi^\#\circ\Phi$ and recall that
$\pi_{AG}^\#\circ j_G' = j_G^{-1}\circ A(\pi^\#)$.
\qed

\begin{thm}                                \label{thm:last}
Let $(\Phi,\phi)$ and $(\Psi,\psi)$ be star 1-forms on $G$, and let
$\omega, \theta\in\fms{P}$. Then
$$
[\Tilde{\Phi},\Tilde{\Psi}] = \Tilde{[\Phi,\Psi]},\qquad
[\Tilde{\Phi},q^*\omega] = q^*(D_\Phi(\omega)),\qquad
[q^*\theta,q^*\omega] = 0.
$$
\end{thm}

\pf
The last equation is of course known. We verify the first two by evaluating
them on all $\Tilde{\zeta}$, where $(\zeta,z)$ is a star vector field on $G$,
and all $Z^\upa$, where $Z\in\Ga AG$. By Proposition \ref{pro:generation},
this suffices.

For any star vector field $(\zeta,z)$, we have
$$
\langle[\Tilde{\Phi},\Tilde{\Psi}],\Tilde{\zeta}\rangle =
\Tilde{\xi}(\langle\Tilde{\Psi},\Tilde{\zeta}\rangle) -
\Tilde{\eta}(\langle\Tilde{\Phi},\Tilde{\zeta}\rangle) +
\langle\Tilde{\Phi},[\Tilde{\eta},\Tilde{\zeta}]\rangle -
\langle\Tilde{\Psi},[\Tilde{\xi},\Tilde{\zeta}]\rangle -
\Tilde{\zeta}(\langle\Tilde{\Phi},\Tilde{\xi}\rangle)
$$
where $\xi = \Phi^\#,\eta = \Psi^\#$. Using (\ref{eq:tilphitilxi}) and
Proposition \ref{pro:sTildes}, it is straightforward to verify that this is
the tilde of $\langle[\Phi,\Psi],\zeta\rangle$, and is therefore equal to
$\langle\Tilde{[\Phi,\Psi]},\Tilde{\zeta}\rangle.$

For any $Z\in\Ga AG$, (\ref{eq:tilphixupa}) gives
$\langle\Tilde{[\Phi,\Psi]},Z^\upa\rangle =
\langle[\phi,\psi],Z\rangle\circ q.$ Similarly expanding out the left hand
side, the equation reduces to
\begin{equation}                          \label{eq:lbabrack}
\langle[\phi,\psi],Z\rangle =
a_*(\phi)(\langle\psi,Z\rangle) - a_*(\psi)(\langle\phi,Z\rangle) +
\langle\phi,D_\eta(Z)\rangle - \langle\psi,D_\xi(Z)\rangle -
Z(\langle\Psi,\xi\rangle).
\end{equation}
Now the equality of these two functions on $P$ follows from similarly
expanding out $\langle[\Phi,\Psi],\Vec{Z}\rangle$ on $G$ and restricting the
result to the identity elements of $G$.

For the second equation, we expand out in the same way, and use the lemma
below.

Lastly, it is easy to see that both sides of
$$
\langle[\Tilde{\Phi},q^*\omega],Z^\upa\rangle =
\langle q^*D_\Phi(\omega),Z^\upa\rangle
$$
are zero, for all $Z\in\Ga AG$.
\qed

\begin{lem}
Let $(\Phi,\phi)$ be a star 1-form, $(\zeta,z)$ a star vector field,
and $\omega\in\fms{P}$. Then
$$
\langle D_\Phi(\omega),z\rangle - \langle\phi,D_\zeta(a_*^*\omega)\rangle
= a_*^*(\omega)\langle\Phi,\zeta\rangle + \delta\omega(a_*\phi,z).
$$
\end{lem}

\pf
Let $X=a_{*}^{*}\omega \in \Ga (AG)$ and let $\Vec{X}$ be the corresponding
right invariant vector field. Thus
\begin{equation}
\label{eq:1}
(a_{*}^{*}\omega )\langle \Phi , \zeta \rangle
= \langle L_{\Vec{X}}\Phi , \zeta\rangle |_{P} +
   \langle \Phi , L_{\Vec{X}}\zeta \rangle |_{P}.
\end{equation}
The second term here is
$$
\langle \Phi , L_{\Vec{X}}\zeta \rangle |_{P}
= -\langle \Phi ,[\zeta , \Vec{X} ]\rangle |_{P}
= -\langle \phi , [\zeta , \Vec{X} ]\smalcirc 1\rangle
= -\langle \phi , D_{\zeta }X\rangle
= -\langle \phi , D_{\zeta } (a_{*}^{*}\omega )\rangle.
$$

On the other hand,
\be
\langle D_{\Phi}\omega , z\rangle &=&\langle [\Phi ,  \beta^{*}\omega ], \zeta \rangle |_{P}\\
&=&\langle L_{\Phi^{\#}}\beta^{*}\omega -L_{(\beta^{*}\omega)^{\#}}\Phi -
\delta (\pi (\Phi , \beta^{*}\omega )), \zeta \rangle |_{P}.
\ee

Clearly, $(\beta^{*}\omega)^{\#}$ is tangent to the $\alpha$-fibers
and is right invariant. In fact,
$(\beta^{*}\omega)^{\#}=-\Vec{(a_{*}^{*}\omega )}=-\Vec{X}$.
Next,
\be
\langle L_{\Phi^{\#}}\beta^{*}\omega , \zeta \rangle |_{P}&=&
\langle \delta i_{\Phi^{\#}}\beta^{*}\omega +i_{\Phi^{\#}}\delta \beta^{*}\omega , \zeta \rangle |_{P}\\
&=&\langle \delta \pi (\Phi , \beta^{*}\omega ), \zeta \rangle |_{P} +
(\delta \beta^{*}\omega) (\Phi^{\#} , \zeta)|_{P}\\
&=&\langle \delta \pi (\Phi , \beta^{*}\omega ), \zeta \rangle |_{P}+
(\delta \omega) (\beta_{*} \Phi^{\#} , \beta_{*}\zeta)|_{P}\\
&=& \langle \delta \pi (\Phi , \beta^{*}\omega ), \zeta \rangle |_{P}
+(\delta \omega)(a_{*}\phi , z).
\ee
Here, in the last step, we used $\beta_{*}\Phi^{\#}=a_{*}\phi $.
This is because on $P$, $\Phi=\phi$ is conormal to $P$ and then
 $\Phi^{\#}=\phi^{\#}=a_{*}\phi$ is tangent to $P$.

Hence,
\begin{equation}
\label{eq:2}
\langle D_{\Phi}\omega , z\rangle =(\delta \omega)(a_{*}\phi , z)
+\langle L_{\Vec{X}}\Phi , \zeta \rangle |_{P}.
\end{equation}

Subtracting Equation (\ref{eq:2}) from (\ref{eq:1}), one gets
$$
a_{*}^{*}(\omega)\langle\Phi,\zeta)-\langle D_{\Phi}\omega,z\rangle
= \langle\Phi,L_{\Vec{X}}\zeta \rangle|_{P} - (\delta \omega )(a_{*}\phi,z)
= -\langle \phi , D_{\zeta}(a_{*}^{*}\omega )\rangle -
         (\delta  \omega )(a_{*}\phi ,z).
$$

This concludes the proof of the lemma.
\qed

By Proposition \ref{pro:genforms}, the three equations in Theorem
\ref{thm:last} completely determine the bracket structure in $\fms{AG}$.
We could, in fact, define the Poisson structure on $AG$ by means of them,
and obtain the Lie algebroid structure of $A^*G$ as its dual. In any case,
equation (\ref{eq:lbabrack}) now gives an explicit expression for the Lie
bialgebroid bracket.

Given a Poisson Lie group $G$ with dual $G^*$ it is elementary that there
is a bijective correspondence between right--invariant vector fields on
$G^*$ and right--invariant 1--forms on $G$. For general Poisson groupoids,
such a result is not possible.

However, suppose that $G\gpd P$ is a Poisson groupoid with dual $G^*\gpd P$
in the sense of \cite{Weinstein:1988}; thus there is an isomorphism of Lie
algebroids $K\colon A(G^*)\to A^*G$ whose dual $K^*\colon AG\to A^*(G^*)$
is also an isomorphism of Lie algebroids. Let $(\xi,x)$ be a multiplicative
vector field on $G$. Then $H_\xi\in{\cal X}(A^*G)$, the dual of $\Tilde{\xi}
\in{\cal X}(AG)$, induces the morphic vector field
$T(K)^{-1}\circ H_\xi\circ K$ on $A(G^*)$. If $G^*$ is $\alpha$--simply
connected, applying Theorem~\ref{thm:intvf} gives a multiplicative vector
field $\eta\in{\cal X}(G^*)$ such that $T(K)\circ\Tilde{\eta} = H_\xi\circ K$.
Reversing this argument gives the following result.

\begin{pro}                \label{pro:last}
Let $G\gpd P$ and $G^*\gpd P$ be $\alpha$--simply connected Poisson
groupoids in duality. Then the above gives a bijective correspondence
between multiplicative vector fields on $G$ and multiplicative vector fields
on $G^*$.
\end{pro}

For example, let $\Gamma\gpd P$ be an $\alpha$--simply connected
symplectic groupoid realizing a simply--connected Poisson manifold $P$.
Then $\Gamma\gpd P$ and $\Bar{P}\times P\gpd P$ are dual Poisson groupoids
\cite{Weinstein:1988}. From \ref{pro:last} it therefore follows that
every vector field on $P$ induces a multiplicative vector field on $\Gamma$,
and every multiplicative vector field on $\Gamma$ arises this way.

\newcommand{\noopsort}[1]{} \newcommand{\singleletter}[1]{#1}

\end{document}